\documentclass[a4paper,12pt]{amsart}
\usepackage{amsmath}
\usepackage{amsthm}
\usepackage{amssymb}
\usepackage{amscd}

\newtheorem{thm}{Theorem}[section]

\newcommand{\proj}{\mathop{\rm Proj}\nolimits}

\DeclareMathOperator{\Ext}{Ext}
\DeclareMathOperator{\Spec}{Spec}

\newcommand{\caltor}{\mathop{{\mathcal T\!or}}\nolimits}

\DeclareMathOperator{\Coker}{Coker}

\DeclareMathOperator{\length}{length}

\DeclareMathOperator{\Supp}{Supp}

\makeatletter

\@addtoreset{equation}{section}
\makeatother

\newtheorem*{acknow}{Acknowledgment}

\title[
Nef vector bundles on a projective space
with $c_1=3$ and $c_2=8$
]{
Nef vector bundles on a projective space
with first Chern class 3 and second Chern class 8
}

\thanks{
This work was partially supported by 
JSPS KAKENHI (C) Grant Number 15K04810.
}

\author{Masahiro Ohno
}

\address{Graduate School of Informatics and Engineering,
The University of Electro-Communications,
Chofu-shi,
Tokyo, 182-8585 Japan
}

\email{masahiro-ohno@uec.ac.jp}



\begin{document}
\begin{abstract}
We describe nef vector bundles on a projective space
with first Chern class three
and second Chern class eight 
over an algebraically closed field of characteristic zero
by giving them a minimal resolution in terms of a full strong exceptional 
collection of line bundles.
\end{abstract}

\maketitle


\section{Introduction}
This paper is a continuation of \cite{Nefofc1=3OnPN}.
Throughout this paper, as in \cite{Nefofc1=3OnPN}, 
we work over an algebraically closed field $K$ of characteristic zero.
Let $\mathcal{E}$ be a nef vector bundle of rank $r$
on a projective space $\mathbb{P}^n$
with first Chern class $c_1$ and second Chern class $c_2$.
In \cite[Theorem 1.1]{Nefofc1=3OnPN},
we classified such $\mathcal{E}$'s in case $c_1=3$ and $c_2< 8$,
and in \cite[Proposition 1.2]{Nefofc1=3OnPN}, 
we also gave an example of such $\mathcal{E}$'s on a projective plane with $c_1=3$ and $c_2=8$.
In this paper, we 
complete the classification of such $\mathcal{E}$'s with $c_1=3$ and $c_2=8$
by 
giving them a minimal resolution
in terms of a full strong exceptional collection of line bundles.
The 
precise statement
is as follows.
\begin{thm}\label{c_1=3c_2=8}
Let $\mathcal{E}$ be as above.
Suppose that $c_1=3$ and that $c_2=8$.
Then $n=2$ and $\mathcal{E}$ fits in 
an 
exact sequence
\[
0\to \mathcal{O}(-2)^{\oplus 2}\to \mathcal{O}^{\oplus r+1}\oplus \mathcal{O}(-1)\to \mathcal{E}\to 0.
\]
\end{thm}
This implies that the example given in \cite[Proposition 1.2]{Nefofc1=3OnPN} is nothing but the unique type
of nef vector bundles with $c_1=3$ and $c_2=8$.

Note that,
for a nef vector bundle $\mathcal{E}$ with $c_1=3$,
the anti-canonical bundle on $\mathbb{P}(\mathcal{E})$ is ample if $n\geq 3$
and nef if $n\geq 2$. 
Moreover, if $n=2$, it is big if and only if $c_2\leq 8$.
So we can say that
we have classified,
except for the case (11) of \cite[Theorem 1.1]{Nefofc1=3OnPN},
weak Fano manifolds of the form $\mathbb{P}(\mathcal{E})$
where $\mathcal{E}$ is a vector bundle on a projective space $\mathbb{P}^n$
under the assumption that $\mathcal{E}$ is nef and $c_1=3$.
Recall here that a projective manifold $M$ is called weak Fano
if its anti-canonical bundle is nef and big,
and that a vector bundle $\mathcal{F}$ is called  a weak Fano bundle
if $\mathbb{P}(\mathcal{F})$ is a weak Fano manifold.
We hope that the theorem above together with \cite[Theorem 1.1]{Nefofc1=3OnPN}
would be useful for some part of the classification of  weak Fano bundles.

This paper is organized as follows.
We first concentrate our attention to  the case $n=2$. In \S~\ref{Set-up for the two-dimensional case},
we recall and summarize results obtained in \cite{Nefofc1=3OnPN}
by taking into account that we only consider nef vector bundles with $c_1=3$ and $c_2=8$.
In \S~\ref{E has O(1) as a subsheaf}, we show that $\mathcal{E}$ does not contain $\mathcal{O}(1)$
as a subsheaf.
In \S~\ref{E does not contain O(1)}, we first observe that $\mathcal{E}$ must fit in the exact sequence
given in \cite[Proposition 1.2]{Nefofc1=3OnPN}
and then show that $\mathcal{E}$ fits in the exact sequence in the theorem above.
Finally, in \S~\ref{nbigger than 2 does not happen},
we show that the case $n\geq 3$ does not happen.

\begin{acknow}
Deep appreciation goes to the referee for his 
careful reading the manuscript and 
invaluable comments.
\end{acknow}

\subsection{Notation and conventions}\label{convention}
Basically we follow the standard notation and terminology in algebraic
geometry. 
For 
a vector bundle $\mathcal{E}$,
$\mathbb{P}(\mathcal{E})$ denotes $\proj S(\mathcal{E})$,
where $S(\mathcal{E})$ denotes the symmetric algebra of $\mathcal{E}$.
For a 
coherent sheaf $\mathcal{F}$ on a smooth projective variety $X$,
we denote by $c_i(\mathcal{F})$ the $i$-th Chern class of $\mathcal{F}$.
For coherent sheaves $\mathcal{F}$ and $\mathcal{G}$ on $X$,
$h^q(\mathcal{F})$ denotes $\dim H^q(\mathcal{F})$.
Finally we refer to \cite{MR2095472} for the definition
and basic properties of nef vector bundles.

\section{Set-up for the two-dimensional case}\label{Set-up for the two-dimensional case}
In the following, let $\mathcal{E}$ be a nef vector bundle on a projective space $\mathbb{P}^n$ with $c_1=3$ and $c_2=8$.
In this section, we assume that $n=2$.
It follows from \cite[(3.10), (3.11) and (3.12)]{Nefofc1=3OnPN} that 
\begin{gather}
h^1(\mathcal{E}(-2))=5,\label{h1E-2onP2}\\
h^0(\mathcal{E}(-1))-h^1(\mathcal{E}(-1))=-2,\label{RRonP2(-1)}\\
h^0(\mathcal{E})=r+1.\label{RRonP2(0)}
\end{gather}
Note here that, for a nef vector bundle $\mathcal{E}'$ in general,
unlike the case of globally generated vector bundles,
an inequality $h^0(\mathcal{E}')\geq r-1$ does not necessarily imply that $\mathcal{E}'$
fits in an exact sequence of the form
\[
0\to \mathcal{O}^{\oplus r-1}\to \mathcal{E}'\to \mathcal{I}_Z\otimes \det\mathcal{E}'\to 0
\]
for some closed subscheme $Z$ of $\mathbb{P}^2$,
where $\mathcal{I}_Z$ denotes the ideal sheaf of $Z$
(see \cite[\S 13]{Nefofc1=3OnPN} for some examples).
Set 
\[e_{0,1}=h^0(\mathcal{E}(-1)).\]
Then 
\[h^1(\mathcal{E}(-1))=e_{0,1}+2\geq 2.\]
It follows from \cite[(3.13)]{Nefofc1=3OnPN} that 
$
5\geq h^1(\mathcal{E}(-1))
$.
Therefore 
\[0\leq e_{0,1}\leq 3.\]

We 
apply to $\mathcal{E}$ the Bondal spectral sequence 
\cite[Theorem 1]{MR3275418}
\begin{equation}\label{BSSonP2}
E_2^{p,q}=\caltor_{-p}^A(\Ext^q(G,\mathcal{E}),G)
\Rightarrow
E^{p+q}=
\begin{cases}
\mathcal{E}& \textrm{if}\quad  p+q= 0\\
0& \textrm{if}\quad  p+q\neq 0.
\end{cases}
\end{equation}
As we have seen in \cite[\S 3.1 and Lemma 5.1]{Nefofc1=3OnPN},
$E_2^{p,q}$ vanishes unless $(p,q)=(-2,1)$, $(-1,1)$ or $(0,0)$,
and $E_2^{-2,1}$ and $E_2^{-1,1}$ fit in an exact sequence of coherent sheaves
\begin{equation}\label{exactseq}
0\to E_2^{-2,1}\to \mathcal{O}(-3)
\xrightarrow{\nu_2}
\Omega_{\mathbb{P}^2}(1)^{\oplus e_{0,1}}
\to E_2^{-1,1}\to k(w)\to 0
\end{equation}
for some point $w$ in $\mathbb{P}^2$,
where $k(w)$ denotes the residue field of $w$.
Note that this exact sequence is a consequence of the vanishing $H^1(\mathcal{E})=0$,
and recall that $H^1(\mathcal{E})$ vanishes by the Kawamata-Viehweg vanishing theorem since $c_2<9$.
Moreover we have the following exact sequences
\begin{gather}
0\to E_2^{-2,1}\to E_2^{0,0}\to E_3^{0,0}\to 0,\label{E_3^{0,0}definition}\\
0\to E_{3}^{0,0}\to \mathcal{E}\to E_2^{-1,1}\to 0,\label{E_2^{-1,1}quotient}\\
0\to \mathcal{O}^{\oplus 3e_{0,1}}\to 
\mathcal{O}(1)^{\oplus e_{0,1}}\oplus \mathcal{O}^{\oplus r+1}
\to E_2^{0,0}\to 0.\label{E_2^00 exact sequence in dim 2}
\end{gather}
We shall divide the proof according to the value of $e_{0,1}$.

\section{The case $n=2$ and $e_{0,1}>0$}\label{E has O(1) as a subsheaf}
Suppose that $n=2$ and $e_{0,1}>0$.
Since $e_{0,1}>0$ and $h^0(\mathcal{E}(-2))=0$ by the argument in \cite[\S 3]{Nefofc1=3OnPN}, we have an exact sequence
\[
0\to \mathcal{O}(1)\to \mathcal{E}\to \mathcal{F}\to 0
\]
where $\mathcal{F}$ is a torsion-free sheaf with $c_1(\mathcal{F})=2$, $c_2(\mathcal{F})=6$
and $h^0(\mathcal{F}(-1))=e_{0,1}-1$.
Denote by $\mathcal{F}^{\vee\vee}$ the double dual of $\mathcal{F}$,
and consider the quotient $\mathcal{Q}$ of the inclusion $\mathcal{F}\subset \mathcal{F}^{\vee\vee}$:
\[
0\to \mathcal{F}\to \mathcal{F}^{\vee\vee}\to \mathcal{Q}\to 0.
\]
The support of $\mathcal{Q}$ has dimension zero,
and 
its length 
is equal to $-c_2(\mathcal{Q})$. 
By \cite[Lemma 12.1]{Nefofc1=3OnPN}, $\mathcal{F}^{\vee\vee}$ is a nef vector bundle of rank $r-1$ 
with $c_1(\mathcal{F}^{\vee\vee})=2$,
$c_2(\mathcal{F}^{\vee\vee})=6+c_2(\mathcal{Q})$
and $h^0(\mathcal{F}^{\vee\vee}(-1))\geq e_{0,1}-1$.

\subsection{The case $e_{0,1}>1$}
Suppose that $e_{0,1}>1$. Then it follows from \cite[Theorem~6.5]{resolution}
that 
$\mathcal{F}^{\vee\vee}$ is isomorphic to either $\mathcal{O}(2)\oplus\mathcal{O}^{\oplus r-2}$
or $\mathcal{O}(1)^{\oplus 2}\oplus\mathcal{O}^{\oplus r-3}$,
or $\mathcal{F}^{\vee\vee}$ fits in an exact sequence
\begin{equation}\label{c_1=2andO(1)sub}
0\to \mathcal{O}(-1)\to \mathcal{O}(1)\oplus \mathcal{O}^{\oplus r-1}\to \mathcal{F}^{\vee\vee}\to 0.
\end{equation}

Suppose that $\mathcal{F}^{\vee\vee}\cong\mathcal{O}(2)\oplus\mathcal{O}^{\oplus r-2}$.
Since $c_2(\mathcal{F}^{\vee\vee})=0$, the length of $\mathcal{Q}$ is $6$.
Let $\mathcal{G}$ be the image of the composite of the inclusion $\mathcal{F}\to \mathcal{O}(2)\oplus\mathcal{O}^{\oplus r-2}$
and the projection $\mathcal{O}(2)\oplus\mathcal{O}^{\oplus r-2}\to \mathcal{O}^{\oplus r-2}$.
Note that the kernel of the surjection $\mathcal{F}\to \mathcal{G}$ is a subsheaf of $\mathcal{O}(2)$.
Hence it can be written as $\mathcal{I}_Z(2)$ where $\mathcal{I}_Z$ is the ideal sheaf of some closed subscheme $Z$ of $\mathbb{P}^2$.
Now we have the following commutative diagram with exact lows and columns
\[
\begin{CD}
 @.     0             @.      0                                       @.     0            @.      \\
@.     @VVV                @VVV                                            @VVV                @. \\
0@>>> \mathcal{I}_Z(2)@>>>\mathcal{O}(2)                              @>>>\mathcal{O}_Z(2)@>>> 0  \\
@.     @VVV                @VVV                                            @VVV                @. \\
0@>>> \mathcal{F}     @>>>\mathcal{O}(2)\oplus\mathcal{O}^{\oplus r-2}@>>>\mathcal{Q  }   @>>> 0  \\
@.     @VVV                @VVV                                            @VVV                @. \\
0@>>> \mathcal{G}     @>>>\mathcal{O}^{\oplus r-2}                    @>>>\mathcal{Q}_1   @>>> 0  \\
@.     @VVV                @VVV                                            @VVV                @. \\
 @.     0             @.      0                                       @.     0            @.        
\end{CD}
\]
where $\mathcal{Q}_1$ is defined by the diagram above.
Since $\mathcal{O}_Z(2)\to \mathcal{Q}$ is injective, we see that $\dim Z\leq 0$,
and thus $\mathcal{O}_Z(2)\cong \mathcal{O}_Z$.
If $\mathcal{Q}_1\neq 0$, then take a line $L$ intersecting with the support of $\mathcal{Q}_1$.
Then the kernel of the surjection $\mathcal{O}_L^{\oplus r-1}\to \mathcal{Q}_1|_L$ has a negative degree line bundle
as a direct summand,
which implies that some negative degree line bundle is a quotient of $\mathcal{G}|_L$, $\mathcal{F}|_L$ and $\mathcal{E}|_L$.
This contradicts that $\mathcal{E}$ is nef. Hence $\mathcal{Q}_1=0$.
Thus $\mathcal{G}\cong \mathcal{O}^{\oplus r-2}$, $\mathcal{O}_Z\cong \mathcal{Q}$, and $\mathcal{O}_Z$ has length $6$.
Since $h^0(\mathcal{G}(-1))=0$, 
we infer that $h^0(\mathcal{I}_Z(1))
=e_{0,1}-1>0$.
Hence there exists a line $L$ passing through $Z$.
Since $\length \mathcal{O}_Z=6$, this implies that
the kernel of the restriction $\mathcal{O}_L(2)\to \mathcal{O}_Z$
to the line $L$ of the surjection $\mathcal{O}(2)\to \mathcal{O}_Z$
is isomorphic to $\mathcal{O}_L(-4)$.
By restricting the diagram above to the line $L$, we see that 
$\mathcal{F}|_L$ has a negative degree line bundle as a quotient;
this is a contradiction.
Hence $\mathcal{F}^{\vee\vee}$ cannot be isomorphic to $\mathcal{O}(2)\oplus\mathcal{O}^{\oplus r-2}$.

Suppose that $\mathcal{F}^{\vee\vee}\cong\mathcal{O}(1)^{\oplus 2}\oplus\mathcal{O}^{\oplus r-2}$.
Since $c_2(\mathcal{F}^{\vee\vee})=1$, the length of $\mathcal{Q}$ is $5$.
Let $\mathcal{G}$ be the image of the composite of the inclusion $\mathcal{F}\to \mathcal{O}(1)^{\oplus 2}\oplus\mathcal{O}^{\oplus r-3}$
and the projection $\mathcal{O}(1)^{\oplus 2}\oplus\mathcal{O}^{\oplus r-3}\to \mathcal{O}^{\oplus r-3}$,
and $\mathcal{Q}_1$ the cokernel of the inclusion $\mathcal{G}\to \mathcal{O}^{\oplus r-3}$.
Then there exists a surjection $\mathcal{Q}\to \mathcal{Q}_1$, and thus the support of $\mathcal{Q}_1$ 
has dimension $\leq 0$.
If $\mathcal{Q}_1\neq 0$, we get a contradiction by the same argument as above.
Therefore we may assume that $\mathcal{Q}_1=0$; thus $\mathcal{G}\cong \mathcal{O}^{\oplus r-3}$.
Let $\mathcal{H}$ be the kernel of the surjection $\mathcal{F}\to \mathcal{O}^{\oplus r-3}$.
Then we have  the following commutative diagram with exact lows and columns.
\[
\begin{CD}
 @.     0                     @.      0     @.                                                            @.     \\
@.     @VVV                       @VVV                                                      @.                @. \\
0@>>> \mathcal{H}             @>>>\mathcal{O}(1)^{\oplus 2}                              @>>>\mathcal{Q}  @>>>0  \\
@.     @VVV                       @VVV                                                      @|                @. \\
0@>>> \mathcal{F}             @>>>\mathcal{O}(1)^{\oplus 2}\oplus\mathcal{O}^{\oplus r-3}@>>>\mathcal{Q}  @>>>0  \\
@.     @VVV                       @VVV                                                      @.                @. \\
 @.   \mathcal{O}^{\oplus r-3}@=  \mathcal{O}^{\oplus r-3}                               @.               @.     \\
@.     @VVV                       @VVV                                                      @.                @. \\
 @.     0                     @.      0     @.                                                              @.       
\end{CD}
\]
Since $h^0(\mathcal{O}^{\oplus r-3}(-1))=0$, 
we infer that $h^0(\mathcal{H}(-1))
=e_{0,1}-1>0$.
Since $\mathcal{H}(-1)$ is a subsheaf of $\mathcal{O}^{\oplus 2}$, this implies that 
$\mathcal{H}(-1)\cong \mathcal{I}_Z\oplus \mathcal{O}$
and $\mathcal{Q}(-1)\cong \mathcal{O}_Z$ for some $0$-dimensional closed subscheme $Z$ of length $5$ in $\mathbb{P}^2$.
Now take a line $L$ that intersect with $Z$ in length $l\geq 2$.
Then the kernel of $\mathcal{O}_L(1)^{\oplus 2}\to \mathcal{O}_{Z\cap L}(1)$ is of the form $\mathcal{O}_L(1-l)\oplus \mathcal{O}_L(1)$.
This implies that $\mathcal{F}|_L$ has a negative degree line bundle as a quotient, which is a contradiction.
Hence $\mathcal{F}^{\vee\vee}$ cannot be isomorphic to $\mathcal{O}(1)^{\oplus 2}\oplus\mathcal{O}^{\oplus r-3}$ either.

Suppose that $\mathcal{F}^{\vee\vee}$ fits in the exact sequence (\ref{c_1=2andO(1)sub}).
Since $c_2(\mathcal{F}^{\vee\vee})=2$, the length of $\mathcal{Q}$ is $4$.
Define a torsion-free sheaf $\mathcal{F}_0$ as a quotient of $\mathcal{F}^{\vee\vee}$
by an injection $\mathcal{O}(1)\to \mathcal{F}^{\vee\vee}$.
Then $\mathcal{F}_0$ fits in an exact sequence
\[
0\to \mathcal{O}(-1)\to \mathcal{O}^{\oplus r-1}\to \mathcal{F}_0\to 0.
\]\
Let $\mathcal{G}$ be the image of the composite of the inclusion $\mathcal{F}\to \mathcal{F}^{\vee\vee}$
and the projection $\mathcal{F}^{\vee\vee}\to \mathcal{F}_0$.
Since $h^0(\mathcal{F}_0(-1))=0$, we see that $h^0(\mathcal{G}(-1))=0$.
Let $\mathcal{H}$ be the kernel of the surjection $\mathcal{F}\to \mathcal{G}$.
Then we have the following commutative diagram with exact lows and columns
\[
\begin{CD}
 @.     0      @.      0     @.                 0           @.     \\
@.     @VVV           @VVV                     @VVV             @. \\
0@>>> \mathcal{H}@>>>\mathcal{O}(1)        @>>>\mathcal{Q}_2@>>>0  \\
@.     @VVV           @VVV                     @VVV             @. \\
0@>>> \mathcal{F}@>>>\mathcal{F}^{\vee\vee}@>>>\mathcal{Q}  @>>>0  \\
@.     @VVV           @VVV                     @VVV             @. \\
0@>>> \mathcal{G}@>>>\mathcal{F}_0         @>>>\mathcal{Q}_1@>>>0  \\
@.     @VVV           @VVV                     @VVV             @. \\
 @.     0      @.      0     @.                 0           @.       
\end{CD}
\]
where $\mathcal{Q}_1$ and  $\mathcal{Q}_2$ are defined by the diagram above.
Since $h^0(\mathcal{G}(-1))=0$,
we see that $h^0(\mathcal{H}(-1))=e_{0,1}-1>0$.
Since $\mathcal{H}(-1)$ is a subsheaf of $\mathcal{O}$, this implies that 
$\mathcal{H}(-1)$ is $\mathcal{O}$ itself;
thus $\mathcal{Q}_2=0$, $\mathcal{Q}\cong\mathcal{Q}_1$ and $\mathcal{Q}_1$ has length $4$.
As we have seen in the proof of \cite[Theorem~6.4]{resolution}, $\mathcal{F}_0$ is locally free
outside at most one point,
and if $\mathcal{F}_0$ is not locally free at a point $z$, then $\mathcal{F}_0$
is isomorphic to $\mathfrak{m}_z(1)\oplus\mathcal{O}^{\oplus r-3}$,
where $\mathfrak{m}_z$ is the ideal sheaf of $z$, since $n=2$.
Suppose that $\mathcal{F}_0$ is not locally free. Then take a line $L$ passing through $z$
and meeting the support of $\mathcal{Q}_1$.
We see that the surjection $\mathcal{F}_0\to \mathcal{Q}_1$ induces a surjection
$\mathcal{O}_L^{\oplus r-2}\to \mathcal{Q}_1|_L$, whose kernel has a negative degree line bundle
as a quotient, and thus so does $\mathcal{G}|_L$, $\mathcal{F}|_L$ and $\mathcal{E}|_L$.
This is a contradiction.
Suppose that $\mathcal{F}_0$ is locally free.
Then take a line $L$ which intersects with $\mathcal{Q}_1$ in length $l\geq 2$.
Since $\mathcal{F}_0|_L\cong \mathcal{O}_L(1)\oplus \mathcal{O}^{\oplus r-3}$,
we see that $\mathcal{G}|_L$ admits a negative degree line bundle as a quotient; this is a contradiction.
Hence $\mathcal{F}^{\vee\vee}$ cannot fit in the exact sequence (\ref{c_1=2andO(1)sub}).

Therefore  we conclude that the case $e_{0,1}>1$ does not happen.

\subsection{The case $e_{0,1}=1$}
Suppose that $e_{0,1}=1$. 
If the morphism $\nu_2$ in (\ref{exactseq}) is zero,
then $E_2^{-1,1}|_L\cong \Omega_{\mathbb{P}^2}(1)|_L\cong \mathcal{O}_L(-1)\oplus \mathcal{O}_L$
for a line $L$ not containing $w$.
By (\ref{E_2^{-1,1}quotient}),
this implies 
that $\mathcal{E}|_L$ has $\mathcal{O}_L(-1)$
as a quotient;
this is a contradiction.
Hence $\nu_2\neq 0$, and thus $E_2^{-2,1}=0$, $E_2^{0,0}\cong E_3^{0,0}$
by (\ref{E_3^{0,0}definition}),
and 
$E_2^{-1,1}$ fits in an exact sequence
\begin{equation}\label{exactseqnew}
0\to \mathcal{O}(-3)
\xrightarrow{\nu_2}
\Omega_{\mathbb{P}^2}(1)
\to E_2^{-1,1}\to k(w)\to 0.
\end{equation}
We see that $E_2^{-1,1}$ 
is a coherent sheaf of rank one.
Since $E_3^{0,0}$ is torsion-free by (\ref{E_2^{-1,1}quotient}),
so is $E_2^{0,0}$, and thus 
$E_2^{0,0}$ has $\mathcal{O}(1)$ as a subsheaf
and consequently is isomorphic to $\mathcal{O}(1)\oplus \mathcal{O}^{\oplus r-2}$
by (\ref{E_2^00 exact sequence in dim 2}).
Hence the exact sequence (\ref{E_2^{-1,1}quotient}) becomes an exact sequence
\[
0\to \mathcal{O}(1)\oplus \mathcal{O}^{\oplus r-2}
\xrightarrow{\varphi} \mathcal{E}
\to E_2^{-1,1}\to 0.
\]
By taking the dual of $\varphi$ and $(r-1)$-th wedge product of the dual,
we obtain a morphism $\wedge^{r-1}\mathcal{E}^{\vee}\to \mathcal{O}(-1)$.
Let $\mathcal{I}_Z(-1)$ be the image of this morphism,
where $\mathcal{I}_Z$ is the ideal sheaf of a closed subscheme $Z$ of $\mathbb{P}^2$
of dimension $\leq 1$.
Note that $Z$ is the degeneracy locus of $\varphi$
and that if we denote by $\psi$ the induced surjection
$\mathcal{E}\cong \wedge^{r-1}\mathcal{E}^{\vee}\otimes \det\mathcal{E}
\to \mathcal{I}_Z(-1)\otimes\det\mathcal{E}\cong \mathcal{I}_Z(2)$
then $\psi\circ \varphi=0$.

Suppose that the degeneracy locus $Z$ of $\varphi$ has codimension $\geq 2$. Then 
$E_2^{-1,1}$ is torsion-free. This implies that $E_2^{-1,1}\cong \mathcal{I}_Z(2)$ and 
that $\mathcal{E}$ fits in an exact sequence
\[
0\to \mathcal{O}(1)\oplus \mathcal{O}^{\oplus r-2}
\xrightarrow{\varphi} \mathcal{E}
\to \mathcal{I}_Z(2)\to 0.
\]
Note that $\length Z=6$.
Since $\mathcal{E}$ is nef,
$\length (Z\cap L)\leq 2$ for any line $L$ in $\mathbb{P}^2$;
let us call this 
the basic property of $Z$.
Let $p$ be any point in $Z$.
We may assume that $Z$ is in 
an affine open subscheme $\Spec K[x,y]$
and that $p=(0,0)$.
The local ring $\mathcal{O}_{Z,p}$ can be written as $A/I$,
where $A=\hat{\mathcal{O}}_{\mathbb{P}^2,p}=K[[x,y]]$
and $I$ the ideal of $Z$ in the local ring $A$.
Observe here that if $\length (A/I)\leq 4$
and thus the support of $Z$ contains another point
$q\neq p$,
then the basic property of $Z$ implies $I\not\subseteq \mathfrak{m}^2$, 
where $\mathfrak{m}$ denotes the maximal ideal of $A$.
Based on this observation,
we can 
deduce from the basic property of $Z$ 
that 
$I\not\subseteq\mathfrak{m}^2$
without any assumption on $\length (A/I)$.
Now that $Z$ is curvilinear,
after changing coordinates $(x,y)$ if necessary, we may assume that 
$I=\langle y-\varphi(x), x^l\rangle$, where $\varphi(x)=a_2x^2+a_3x^3+\dots \in K[[x]]$ $(a_2\neq 0)$ and $l=\length (A/I)$.
Local computation then shows that there exists a smooth conic $C$ such that $\length (Z\cap C)\geq 5$;
e.g., if $l\geq 3$, we can take a defining equation of $C$ to be  $y=a_2x^2+dxy+ey^2$ for some $d, e\in K$.
However this again contradicts that $\mathcal{E}$ is nef. 
Therefore this case cannot happen.

Suppose that $\dim Z=1$.
Then the ideal sheaf $\mathcal{I}_Z$ of $Z$ is decomposed as $\mathcal{I}_Z\cong \mathcal{I}_{Z_d}(-d)$,
where $d$ is the degree of the divisor contained in $Z$
and $\mathcal{I}_{Z_d}$ is the ideal sheaf of 
a $0$-dimensional closed subscheme $Z_d$ of $\mathbb{P}^2$.
Consider the following commutative diagram with exact lows and columns
\[
\begin{CD}
 @.                                               @.   0                    @.   0                    @.      \\
@.     @.                                             @VVV                      @VVV                       @. \\
0@>>> \mathcal{O}(1)\oplus\mathcal{O}^{\oplus r-2}@>>>\mathcal{K}           @>>>\mathcal{T}           @>>> 0  \\
@.     @|                                             @VVV                      @VVV                       @. \\
0@>>> \mathcal{O}(1)\oplus\mathcal{O}^{\oplus r-2}@>{\varphi}>>\mathcal{E}           @>>>E_2^{-1,1}            @>>> 0  \\
@.     @.                                             @VVV                      @VVV                       @. \\
 @.                                               @.  \mathcal{I}_{Z_d}(2-d)@=  \mathcal{I}_{Z_d}(2-d)@.      \\
@.     @.                                             @VVV                      @VVV                       @. \\
 @.                                               @.   0                    @.     0                  @.        
\end{CD}
\]
where $\mathcal{K}$ and $\mathcal{T}$ are defined by the diagram above.
We see that $\mathcal{K}$ is a coherent sheaf of rank $r-1$
and thus $\mathcal{T}$ is the torsion subsheaf of $E_2^{-1,1}$,
and that $\Supp Z= \Supp \mathcal{T}\cup \Supp Z_d$.
Hence $E_2^{-1,1}$ has an associated point of codimension one.
Now recall the exact sequence~(\ref{exactseqnew})
and split this sequence into the following two exact sequences of coherent sheaves
\begin{gather}
0\to \mathcal{O}(-3)
\xrightarrow{\nu_2}
\Omega_{\mathbb{P}^2}(1)
\to \mathcal{C}\to  0,\\
0\to \mathcal{C}
\to E_2^{-1,1}\to k(w)\to 0.
\end{gather}
Note that $\mathcal{C}$ has an associated point of codimension one
since so does $E_2^{-1,1}$. Hence $\nu_2$ 
passes through $\mathcal{O}(-1)$ or $\mathcal{O}(-2)$.

Suppose that $\nu_2$ passes through $\mathcal{O}(-1)$.
Then we have the following commutative diagram with exact lows and columns
\[
\begin{CD}
 @.    0               @.   0                      @.               @.      \\
@.     @VVV                @VVV                        @.                @. \\
 @.   \mathcal{O}(-3)  @=  \mathcal{O}(-3)         @.               @.      \\
@.     @VVV                @VVV                        @.                @. \\
0@>>> \mathcal{O}(-1)  @>>>\Omega_{\mathbb{P}^2}(1)@>>>\mathcal{I}_p@>>> 0  \\
@.     @VVV                @VVV                        @|                @. \\
0@>>> \mathcal{O}_D(-1)@>>>\mathcal{C}             @>>>\mathcal{I}_p@>>> 0  \\
@.     @VVV                @VVV                        @.                @. \\
 @.    0               @.   0                      @.               @.        
\end{CD}
\]
where $\mathcal{I}_p$ is the ideal sheaf of a point $p$,
and $D$ is a conic in $\mathbb{P}^2$.
We also have the following commutative diagram with exact lows and columns
\[
\begin{CD}
 @.    0               @.   0                      @.               @.      \\
@.     @VVV                @VVV                        @.                @. \\
 @.   \mathcal{O}_D(-1)@=  \mathcal{O}_D(-1)       @.               @.      \\
@.     @VVV                @VVV                        @.                @. \\
0@>>> \mathcal{C}      @>>>E_2^{-1,1}              @>>>k(w)         @>>> 0  \\
@.     @VVV                @VVV                        @|                @. \\
0@>>> \mathcal{I}_p    @>>>\mathcal{D}             @>>>k(w)         @>>> 0  \\
@.     @VVV                @VVV                        @.                @. \\
 @.    0               @.   0                      @.               @.        
\end{CD}
\]
where $\mathcal{D}$ is defined by the diagram above.
Suppose that $\mathcal{D}$ has an associated point other than the generic point.
Then it must be $w$, and thus $\mathcal{D}\cong \mathcal{I}_w\oplus k(w)$,
which also contradicts that $\mathcal{E}$ is nef. Therefore $\mathcal{D}$ is torsion-free.
Since $\mathcal{D}$ has rank one, $c_1(\mathcal{D})=0$ and $c_2(\mathcal{D})=0$,
$\mathcal{D}$ is isomorphic to its double dual $\mathcal{O}_{\mathbb{P}^2}$.
Moreover we see that $p=w$,
that $\mathcal{O}_D(-1)$ is the torsion subsheaf $\mathcal{T}$ of $E_2^{-1,1}$,
that $Z_d=\emptyset$, and that $Z=D$.
If $h^0(E_2^{-1,1})\neq 0$, then $E_2^{-1,1}\cong \mathcal{O}_D(-1)\oplus \mathcal{O}_{\mathbb{P}^2}$,
which contradicts that $\mathcal{E}$ is nef. 
Hence $h^0(E_2^{-1,1})= 0$.
Since $h^1(\mathcal{O}_D(-1))=h^2(\mathcal{O}_{\mathbb{P}^2}(-3))=1$,
this implies that $H^0(\mathcal{D})=H^0(\mathcal{O}_{\mathbb{P}^2})\cong H^1(\mathcal{O}_D(-1))$.
Suppose that $D$ is smooth.
Consider the pull back 
$
\mathcal{O}_D(-1)\to E_2^{-1,1}|_D\to \mathcal{O}_D\to 0
$
of the exact sequence above.
Note that $E_2^{-1,1}|_D$ has rank at least two
since $D$ is the degeneracy locus of $\varphi$.
Hence we obtain an exact sequence
\[
0\to \mathcal{O}_D(-1)\to E_2^{-1,1}|_D\to \mathcal{O}_D\to 0.
\]
Note that $D\cong\mathbb{P}^1$
and that $\mathcal{O}_D(-1)\cong \mathcal{O}_{\mathbb{P}^1}(-2)$
via this isomorphism.
Since the sequence above does not split,
this implies that $E_2^{-1,1}|_D\cong \mathcal{O}_{\mathbb{P}^1}(-1)^{\oplus 2}$,
which contradicts that $\mathcal{E}$ is nef.
Suppose that $D$ is a double line.
Then $D_{\textrm{red}}\cong \mathbb{P}^1$,
and we have a surjection 
$\mathcal{O}_D(-1)\to \mathcal{O}_{D_{\textrm{red}}}(-1)$.
The similar argument as above shows that there exists an exact sequence
\[
0\to \mathcal{O}_{D_{\textrm{red}}}(-1)\to E_2^{-1,1}|_{D_{\textrm{red}}}\to \mathcal{O}_{D_{\textrm{red}}}\to 0.
\]
Hence $E_2^{-1,1}|_{D_{\textrm{red}}}\cong \mathcal{O}_{\mathbb{P}^1}(-1)\oplus\mathcal{O}_{\mathbb{P}^1}$;
this contradicts that $\mathcal{E}$ is nef.
Suppose that $D$ is a union of two distinct lines: $D=L_1+L_2$.
Then 
$E_2^{-1,1}|_{L_1}\cong \mathcal{O}_{\mathbb{P}^1}(-1)\oplus\mathcal{O}_{\mathbb{P}^1}$
by the similar argument as above,
and hence this case does not occur either.

Suppose that $\nu_2$ passes through $\mathcal{O}(-2)$
and does not pass through $\mathcal{O}(-1)$.
Then we have the following commutative diagram with exact lows and columns
\[
\begin{CD}
 @.    0               @.   0                      @.                  @.      \\
@.     @VVV                @VVV                        @.                   @. \\
 @.   \mathcal{O}(-3)  @=  \mathcal{O}(-3)         @.                  @.      \\
@.     @VVV                @VVV                        @.                   @. \\
0@>>> \mathcal{O}(-2)  @>>>\Omega_{\mathbb{P}^2}(1)@>>>\mathcal{I}_W(1)@>>> 0  \\
@.     @VVV                @VVV                        @|                   @. \\
0@>>> \mathcal{O}_L(-2)@>>>\mathcal{C}             @>>>\mathcal{I}_W(1)@>>> 0  \\
@.     @VVV                @VVV                        @.                   @. \\
 @.    0               @.   0                      @.                  @.        
\end{CD}
\]
where $\mathcal{I}_W$ is the ideal sheaf of a $0$-dimensional locally complete intersection 
$W$ of length three,
and $L$ is a line in $\mathbb{P}^2$.
We also have the following commutative diagram with exact lows and columns
\[
\begin{CD}
 @.    0               @.   0                      @.               @.      \\
@.     @VVV                @VVV                        @.                @. \\
 @.   \mathcal{O}_L(-2)@=  \mathcal{O}_L(-2)       @.               @.      \\
@.     @VVV                @VVV                        @.                @. \\
0@>>> \mathcal{C}      @>>>E_2^{-1,1}              @>>>k(w)         @>>> 0  \\
@.     @VVV                @VVV                        @|                @. \\
0@>>> \mathcal{I}_W(1) @>>>\mathcal{D}             @>>>k(w)         @>>> 0  \\
@.     @VVV                @VVV                        @.                @. \\
 @.    0               @.   0                      @.               @.        
\end{CD}
\]
where $\mathcal{D}$ is defined by the diagram above.
If $\mathcal{D}$ is not torsion-free, then $\mathcal{D}\cong \mathcal{I}_w\oplus k(w)$,
which contradicts that $\mathcal{E}$ is nef. 
Therefore $\mathcal{D}$ is a torsion-free coherent sheaf of rank one with
$c_1(\mathcal{D})=1$ and $c_2(\mathcal{D})=2$.
Hence $\mathcal{O}_L(-2)$ is the torsion subsheaf $\mathcal{T}$ of $E_2^{-1,1}$,
and we infer that $\mathcal{D}\cong \mathcal{I}_{Z_1}(1)$ with $\length Z_1=2$.
This also contradicts that $\mathcal{E}$ is nef.

Therefore we conclude that the case $e_{0,1}=1$ does not happen.

\section{The case $n=2$ and $e_{0,1}=0$}\label{E does not contain O(1)}
Suppose that $e_{0,1}=0$.
Then $E_2^{-2,1}\cong \mathcal{O}(-3)$ and $E_2^{-1,1}\cong k(w)$ by (\ref{exactseq}),
and $E_2^{0,0}\cong \mathcal{O}^{\oplus r+1}$ by (\ref{E_2^00 exact sequence in dim 2}).
Thus we have the following two exact sequences by (\ref{E_3^{0,0}definition}) and (\ref{E_2^{-1,1}quotient})
\begin{gather}
0\to \mathcal{O}(-3)\to \mathcal{O}^{\oplus r+1}\to E_3^{0,0}\to 0,\\
0\to E_{3}^{0,0}\to \mathcal{E}\to k(w)\to 0.
\end{gather}
These two exact sequences show that $\mathcal{E}$ must fit in the exact sequence
given in \cite[Proposition 1.2]{Nefofc1=3OnPN}.
We shall show that $\mathcal{E}$ has a resolution in terms of a full strong exceptional sequence 
of line bundles as in Theorem~\ref{c_1=3c_2=8} in accordance with the framework given in \cite{resolution}.

Since $h^1(E_3^{0,0}(1))=0$, we have the following commutative diagram with exact rows and columns
\[
\begin{CD}
 @.    0                     @.   0                                           @.   0               @.      \\
@.    @VVV                       @VVV                                             @VVV                  @. \\
0@>>>\mathcal{O}(-3)         @>>>\mathcal{J}                                  @>>>\mathcal{I}_w(-1)@>>> 0  \\
@.    @VVV                       @V{g}VV                                             @VVV                  @. \\
0@>>>\mathcal{O}^{\oplus r+1}@>>>\mathcal{O}^{\oplus r+1}\oplus\mathcal{O}(-1)@>>>\mathcal{O}(-1)  @>>> 0  \\
@.    @VVV                       @VVV                                             @VVV                  @. \\
0@>>>E_3^{0,0}               @>>>\mathcal{E}                                  @>>>k(w)             @>>> 0  \\
@.    @VVV                       @VVV                                             @VVV                  @. \\
 @.    0                     @.   0                                           @.   0               @.        
\end{CD}
\]
where $\mathcal{I}_w$ is the ideal sheaf of $w$, and $\mathcal{J}$ and $g$ are defined by the diagram above.
We also have the following commutative diagram with exact rows and columns
\[
\begin{CD}
 @.                 @.  0                                              @.  0                         @.      \\
@.   @.                 @VVV                                               @VVV                           @. \\
 @.                 @.  \mathcal{O}(-3)                                @=  \mathcal{O}(-3)           @.      \\
@.   @.                 @V{f}VV                                               @VVV                           @. \\
0@>>>\mathcal{O}(-3)@>>>\mathcal{O}(-3)\oplus\mathcal{O}(-2)^{\oplus 2}@>>>\mathcal{O}(-2)^{\oplus 2}@>>> 0  \\
@.   @|                 @VVV                                               @VVV                           @. \\
0@>>>\mathcal{O}(-3)@>>>\mathcal{J}                                    @>>>\mathcal{I}_w(-1)         @>>> 0  \\
@.   @.                 @VVV                                               @VVV                           @. \\
 @.                 @.  0                                              @.  0                         @.        
\end{CD}
\]
where $f$ is defined by the diagram above.
We claim here that the composite of $f$ and the projection 
$\mathcal{O}(-3)\oplus\mathcal{O}(-2)^{\oplus 2}\to \mathcal{O}(-3)$ is non-zero.
Suppose, to the contrary, that the composite is zero.
Then $\mathcal{J}\cong \mathcal{O}(-3)\oplus \mathcal{I}_w(-1)$.
By taking the double dual, the composite of the inclusion $\mathcal{I}_w(-1)\to \mathcal{J}$ and $g$
extends to a splitting injection of the projection $\mathcal{O}^{\oplus r+1}\oplus\mathcal{O}(-1)\to \mathcal{O}(-1)$;
we obtain the following commutative diagram with exact rows
\[
\begin{CD}
0@>>>\mathcal{I}_w(-1)  @>>>\mathcal{O}(-1)                                 @>>>k(w)       @>>> 0  \\
@.   @VVV               @VVV                                                @VVV            @. \\
0@>>>\mathcal{J}    @>{g}>>\mathcal{O}^{\oplus r+1}\oplus\mathcal{O}(-1)@>>>\mathcal{E}@>>> 0.  \\
\end{CD}
\]
Since the induced morphism $k(w)\to \mathcal{E}$ is a splitting injection of the surjection $\mathcal{E}\to k(w)$,
we have an isomorphism $\mathcal{E}\cong E_3^{0,0}\oplus k(w)$, which is absurd.
Hence the claim holds; thus 
$\mathcal{J}\cong \Coker(f)\cong \mathcal{O}(-2)^{\oplus 2}$.
Therefore we obtain the desired exact sequence
\[
0\to \mathcal{O}(-2)^{\oplus 2}\to \mathcal{O}^{\oplus r+1}\oplus \mathcal{O}(-1)\to \mathcal{E}\to 0.
\]

\section{The case $n\geq 3$}\label{nbigger than 2 does not happen}
In this section, we shall show that the case $n\geq 3$ does not happen.
By considering the restriction $\mathcal{E}|_{L^3}$ to a $3$-dimensional 
linear subspace $L^3\subseteq \mathbb{P}^n$,
we may assume that $n=3$.
We have 
\[
\chi(\mathcal{E}(-1))
=\frac{c_3}{2}-10
\]
by \cite[(3.20)]{Nefofc1=3OnPN}.
In particular, $c_3$ is even.
We also have 
\[
c_3\geq 
21
\]
by \cite[(3.23)]{Nefofc1=3OnPN}.
Since the equality in $c_3\geq 21$ does not hold,
we infer that $H(\mathcal{E})$ is big, and thus $h^q(\mathcal{E}(-1))=0$ for all $q>0$
by \cite[(3.3)]{Nefofc1=3OnPN}.
Therefore $h^0(\mathcal{E}(-1))\geq 1$.
On the other hand, 
$H^0(\mathcal{E}(-2))=0$ by 
the argument in \cite[\S 3]{Nefofc1=3OnPN},
and $h^0(\mathcal{E}|_H(-1))=0$ for any plane $H\subset \mathbb{P}^3$
as is shown in \S~\ref{E has O(1) as a subsheaf}.
Hence $h^0(\mathcal{E}(-1))=0$, which is a contradiction.
Therefore the case $n\geq 3$ does not happen.

\bibliographystyle{alpha}
\bibliography{c2=8NefRevised.bbl}
\end{document}